\documentclass[12pt,leqno]{article}

\usepackage{amsmath,amssymb}
\usepackage{latexsym}
\usepackage{graphicx}
\usepackage{bm}

\def\R{\mathbb{R}}
\def\C{\mathbb{C}}

\begin{document}

\baselineskip=16pt

\title{The Pompeiu problem}

\author{ A. G. Ramm
}

\renewcommand{\thefootnote}{\fnsymbol{footnote}}

\date{}
\maketitle
\begin{abstract}
\baselineskip=14pt
Let $f \in L_{loc}^1 (\R^n)\cap \mathcal{S}'$, where
$\mathcal{S}'$ is the Schwartz class of distributions,  and
$$\int_{\sigma (D)} f(x) dx = 0 \quad \forall \sigma \in G, \qquad (*)$$
where  $D\subset \R^n$ is a bounded domain, the closure $\bar{D}$ of
which is diffeomorphic to a closed ball. Then the complement of $\bar{D}$
is connected and path connected.
Here  $G$ denotes the group of all rigid motions in $\R^n$.
This group
consists of all translations and rotations.

It is conjectured that if $f\neq 0$ and (*) holds, then $D$ is a ball.
Two other conjectures, equivalent to the above one, are formulated
and discussed.

\end{abstract}

\textbf{MSC:} 35J05, 31B20

\textbf{Key words:} Symmetry problems; Pompeiu problem.

\section{Introduction}
In this paper the problem known as the Pompeiu problem is formulated and
discussed. This problem originated in Pompeiu's paper \cite{P}, of 1929.
The problem in a modern formulation is stated below as Conjecture 1, and
is still open.

Dimitrie  Pompeiu (1873-1954) was born in Romania and got his Ph.D in 1905 
at
the Sorbonne, in Paris, under the direction of H.Poincar$\acute{e}$. He is 
known
mainly for the Pompeiu problem and for the Cauchy-Pompeiu formula in
complex analysis.

Let us formulate the Pompeiu problem as it is understood today.

Let $f \in L_{loc}^1 (\R^n)\cap \mathcal {S}'$, where
$\mathcal{S}'$ is the Schwartz class of distributions,  and
\begin{equation}
\label{e1}
\int_{\sigma (D)} f(x) dx = 0 \quad \forall \sigma \in G,
\end{equation}
where $G$  is the group of all rigid motions of $\R^n$, consisting of
all translations and rotations, and
 $D\subset \R^n$ is a bounded domain, the closure $\bar{D}$ of
which is diffeomorphic to a closed ball. Under these assumptions the
complement of $\bar{D}$
in $\R^n$ is connected and path connected by the isotopy
extension theorem, see \cite{H}.

The distribution space
$\mathcal{S}'$ in the
assumption
$f \in L_{loc}^1 (\R^n)\cap \mathcal{S}'$ can be replaced without
essential changes in the argument by the distribution space
$\mathcal{D}'$, where $\mathcal{D}$ is the space of $C^\infty_0(\R^3)$
functions.

In \cite{P} the following question was raised:

{\bf Does (1) imply that  $f = 0$?}

If yes, then we say that  {\it $D$  has $P$-property} (Pompeiu's
property), and write $D\in P$. Otherwise, we say
that {\it $D$ fails to have $P$-property}, and write  $D\in \overline{P}$.
Pompeiu claimed in 1929 that every plane bounded domain has 
$P$-property.
This claim turned out to be false:
a counterexample was given 15 years later in \cite{C}. The counterexample
is a
domain $D$ which is a disc, or a ball in $\R^n$ for $n> 2$. If $D$ is a
ball, then there are $f\neq 0$ for which equation (1) holds.
The set of all $f\neq 0$, for which equation
(1) holds, was constructed in \cite{R363}. There are infinitely many
(a continuum)  such $f$. Let us give the counterexample mentioned above.

{\bf Example 1.} {\it Suppose that $D\subset \R^n$ is a ball $B$ 
centered at
the origin and of
radius $a$. Then the Fourier transform of its characteristic function
$\chi$ is
\begin{equation}
\label{e0}
\tilde{\chi}(\xi)=\int_Be^{i\xi\cdot
x}dx=(2\pi a)^{n/2}\frac{J_{n/2}(a|\xi|)}
{|\xi|^{n/2}},
\end{equation}
where $J_{n/2}$ is the Bessel function.
It follows that if $|\xi|=s_{j,n}/a$, where $s_{j,n}$ is any positive
zero
of the Bessel function $J_{n/2}(s)$, then $\tilde{\chi}(\xi)$
has a spherical set of zeros. This implies, as follows from Theorem 3,
proved below, that there are $f\neq 0$ for which relation
\eqref{e1} holds.
}
\vskip .2 in
A bibliography on the Pompeiu
problem ($P$-problem) can be found in \cite{Z} and
in \cite{R363}.

The current formulation of the $P$-problem is the following:

{\it Prove that if $D\subset \R^n$ is a bounded domain diffeomorphic
to a ball and $D\in \overline{P}$,
then $D$ is a ball.}

We use the word ball also in the case $n=2$, when this word means
disc, and discuss the $P$-problem in detail. This problem leads
to some problems of general mathematical interest: a symmetry problem
for partial differential equations, see Conjecture 2 below, and a problem 
in harmonic analysis, see Conjecture 3 below.

Let us make the following standing assumptions:

{\it Assumptions A:

$A_1):$ $D$ is a bounded domain, the closure of which is diffeomorphic to
a closed ball, the boundary $S$ of $D$ is a closed connected $C^1$-smooth
surface,

$A_2):$ $D$ fails to have $P$-property.
}

Our first conjecture is:

{\bf Conjecture 1.} {\it If Assumptions A hold, then $D$ is a ball.}

In Section 2 this Conjecture  is discussed. We prove that  Conjecture 1
is equivalent to a symmetry problem for a partial differential equation.
Namely, it is equivalent to the following conjecture.

{\bf Conjecture 2.} {\it If problem \eqref{e2} (see below) has a solution,
then $D$ is a ball.}

Several symmetry problems were studied by the method similar
to the one used in the proof of Theorem 1, below,
see also \cite{R512}-\cite{R614}.

Conjectures 1 and 2 are  equivalent to the following conjecture:

{\bf Conjecture 3.} {\it  If Assumption $A_1$ holds and the Fourier
transform $\tilde{\chi}_D$ of the characteristic function $\chi_D$ of 
the domain $D$ has a spherical surface of
zeros, then $D$ is a ball.}

\section{Discussion of the Conjectures}

It is proved in \cite{W} that if Assumptions A hold, then the boundary
$S$ of $D$ is real-analytic. It is proved in Theorem 3 below,  that
if Assumptions A hold, then  the problem
\begin{equation}
\label{e2}
 (\nabla^2 + k^2)u = 1\quad \text{in}\quad D,\qquad u\big{|}_S = 0,
\quad u_N\big{|}_S=0, \quad k^2 = const>0,
\end{equation}
has a solution.
In \eqref{e2} $N=N_s$ is the outer unit normal to $S$
pointing out
of $D$, $s\in S$ is a point on $S$.

Therefore, if \eqref{e1} holds, then problem \eqref{e2} has a solution.

Let us prove that if  problem \eqref{e2} has a solution, then 
$\tilde{\chi}_D$ has a
spherical set of zeros, where $\chi_D$ is the characteristic function of
the domain $D$. To prove this, 
let us multiply \eqref{e2} by $e^{ik\alpha \cdot x},$ where $\alpha\in 
S^2$ is an arbitrary unit vector and $S^2$ is the unit sphere in $\R^3$,
and integrate with respect to $x$ over $D$. Using an integration by parts 
and the boundary conditions  \eqref{e3} for $u$, one gets the desired
relation:
\begin{equation}
\label{e2a}
\tilde{\chi}_D(k\alpha)=0 \quad  \forall \alpha\in S^2.
\end{equation}
Thus,  $\tilde{\chi}_D$  has a spherical set of zeros.

Conjectures 1, 2, and 3 are equivalent in the sense that each of them 
claims that $D$ is a ball.
But these three Conjectures are also equivalent in the sense that if
one of them is correct, then so are the remaining two.

Indeed, in the proof of Theorem 3 equation \eqref{e2} was derived
from equation \eqref{e1}, and the relation
\begin{equation} \label{eC3}
\tilde{\chi}_D(\xi)=0, \qquad |\xi|=const>0,
\end{equation}
 was derived from equation \eqref{e2}. From the proof of Theorem 3
it follows that  equation \eqref{eC3} implies the equation
\begin{equation} \label{eC4}
\tilde{\chi}(\xi)=(\xi^2-k^2)\tilde{u}(\xi),
\end{equation}
(see  equation \eqref{eZ} below, in the proof of Theorem 3), and the 
inverse
Fourier transform of \eqref{eC4} yields equation \eqref{e2} and
the boundary conditions in \eqref{e2}. In this sense problem
\eqref{e2} is equivalent to relation \eqref{eC3}, and Conjectures 2 and 3 
are equivalent in the sense that if one of them is correct, so is the
other one. 

Furthermore, from the proof of Theorem 3 it follows that
\eqref{eC3} implies \eqref{e1}. Indeed, equation \eqref{eC3} implies
equation 
\begin{equation} \label{eC5}
\tilde{f}(\xi)\overline{\tilde{\chi}(g^{-1}\xi)}=0,
\end{equation}
(see equation \eqref{eX} below), and the inverse Fourier transform of this 
equation yields relation \eqref{e1}. 

In this sense relations \eqref{e1}, \eqref{e2} and \eqref{eC3} are
equivalent, and Conjectures 1, 2 and 3, are equivalent in the sense
that if one of them is correct, so are the other two.

Problem \eqref{e2} is an open symmetry problem of long standing for 
partial differential equations. 
Let us formulate another open symmetry problem for partial differential 
equations of long standing, known as
M. Schiffer's conjecture:

{\bf Conjecture 4.} {\it If the problem
\begin{equation}
\label{eS}
(\nabla^2 + k^2)u = 0\quad \text{in}\quad D,\qquad u_N\big{|}_S = const\neq 0,
\quad u\big{|}_S=0, \quad k^2 = const>0,
\end{equation}
has a solution, then $D$ is a ball.}

Note that Conjecture 2 can be formulated in the form similar to 
\eqref{eS}:

{\bf Conjecture 5.} {\it If the problem
\begin{equation}
\label{eP}
(\nabla^2 + k^2)u = 0\quad \text{in}\quad D,\qquad u_N\big{|}_S = 0,
\quad u\big{|}_S=const\neq 0, \quad k^2 = const>0,
\end{equation}
has a solution, then $D$ is a ball.}

Conjecture 5 is equivalent to Conjecture 2. Indeed, if \eqref{e2}
holds, then
one can look for the solution $u$ of the form $u=v+c$, where $c$ is a 
constant. The boundary
conditions in \eqref{e2} imply $v=-c$ on
$S$, and  $v_N=0$ on $S$. Let us choose $c=1/k^2$. Then the 
differential 
equation 
\eqref{e2} implies  
$$(\nabla^2 + k^2)v=1-k^2c=0.$$
Therefore,  $v$ solves problem \eqref{eP}. Conversely, if $v$ solves 
problem
\eqref{eP}, then $u=v+c$  solves 
\eqref{e2} if $c$ is a suitable constant. \hfill $\Box$

Conjecture 5 is {\it not} equivalent to Conjecture 4.

 The results, on which our discussion of
Conjectures 1, 2 and 3 is based, are Theorems 1, 2 and 3. Theorems 1
and 2 were proved originally in \cite{R382}, \cite{R363},  and in
the book \cite{R470}, Chapter 11. A result, equivalent to Theorem 3, 
had 
been proved originally in the paper \cite{BST} by a considerably longer 
and more complicated argument. Our proofs are borrowed essentially from
\cite{R470}, Chapter 11. In the paper \cite{K} the null-varieties of
the Fourier transform of the characteristic function of a bounded
domain $D$  are studied. The properties of these varieties and the
geometrical properties of $D$ are related, of course, but it is not
clear in what way they are related. Conjecture 3, if it is  proved, is an 
interesting example of such a relation.

In Section 3 a relation of the Pompeiu problem in $\R^2$ to 
analyticity of $f$ is discussed. It is proved that if 
the domain $D\in P$, $f\in L^1_{loc}(\R^2)$, and if $\int_{\partial 
\sigma(D)}fdz=0 \quad \forall \sigma\in G$, then $f$ is an entire 
function. An earlier discussion of this result can be found in \cite{Z1},
and a new short proof of a result from \cite{Z1} is given.

In Section 4 a new approach to the Pompeiu problem is outlined
and a new Conjecture is formulated.  

To make
our presentation essentially self-contained,  proofs of
Theorems 1, 2 and 3,  are included in this paper.

{\bf Theorem 1.}  {\it If Assumptions A hold, then
\begin{equation}
\label{e3}
[s,N]=u_N, \quad \forall s\in S,
\end{equation}
where $[s,N]$ is the cross product in $\R^3$, and $u$ is a
vector-function that solves the problem
\begin{equation}
\label{e4}
(\nabla^2+k^2)u=0 \quad in \quad D, \quad u|_S=0.
\end{equation}
}
If $n=2$, then $D$ is a plane domain, $S$ is a curve, diffeomorphic
to a circle, $u$ is a scalar solution to equation (4), and equation (3)
yields $s_1N_2-s_2N_1=u_N, \quad \forall s\in S,$
where $N_j$, $j=1,2,$ are Cartesian coordinates of the unit normal
$N$ to $S$.
Indeed, if $n=2$ then the cross product of two vectors
$[s_1e_1+s_2e_2, N_1e_1+N_2e_2]$ is calculated by the formula
$[s,N]=(s_1N_2-s_2N_1)e_3,$ where $e_3$ is a unit vector,
orthogonal to the plane domain $D$, and the triple $\{e_j\}_{j=1}^3$
is a standard orthonormal basis in $\R^3$.

Let us state the following characterization of spheres.

{\bf Theorem 2.} {\it If $S$ is a smooth surface homeomorphic to a sphere
and $[s,N]=0$ on $S$, then $S$ is a sphere.}

The proof of Theorem 2 will be given in the coordinate system in which
the condition $[s,N]=0$ on $S$ is valid. 


The following conclusion is an immediate consequence of
Theorem 2:

{\it The conclusion of Conjecture 1 will
be established if one proves, under the Assumptions A,
 that  $[s,N]=0$ on $S$}.

Let us start by  proving Theorem 2, then Theorem 1 is proved, and,
finally, we prove Theorem 3.
It is assumed throughout, except in Section 3, that $n=3$. Our 
proofs of Theorems 2 and 3 can be used
for any integer $n\ge 2$ without any essential changes.
The proof of Theorem 1 uses the notion of the cross-product,
and by this reason its proof should be modified for $n>3$.

{\it Proof of Theorem 2.} Let $n=3$ and
assume that $s=s(p,q)$ is a parametric equation of the surface $S$.
The normal $N$ to $S$ is a vector directed along the vector
$[s_{p},s_{q}]$,
where $s_{p}$ denotes the partial derivative of the function $s(p,q)$ with 
respect to the
parameter $p$.  The assumption $[s, N]=0$ on $S$, yields
\begin{equation}
\label{e5}
[s, [s_p, s_q]]=s_p s\cdot s_q-s_q s\cdot s_p=0,
\end{equation}
where $ s\cdot s_q$ is the dot product of two vectors in $\R^3$.
At a non-singular points of $S$, the vectors $s_p$ and $s_q$ are 
linearly independent. The surface $S$ is smooth, so its points are 
non-singular.  Therefore
equation \eqref{e5}
implies $ s\cdot s_q=0$ and $ s\cdot s_p=0$, so
\begin{equation}
\label{e6}
\frac{\partial s\cdot s}{\partial p}=0,\quad \frac{\partial s\cdot s}
{\partial q}=0.
\end{equation}
Therefore
\begin{equation}
\label{e7}
s\cdot s=const.
\end{equation}
This is an equation of a sphere in the coordinate system with the origin
at the center of the sphere.
Theorem 2 is proved.\hfill $\Box$

{\it Proof of Theorem 1.} Let $n=3$ and $\mathcal{N}$ denote the set of
all smooth solutions to \eqref{e4} in a ball $B$, containing $D$.
Multiply \eqref{e2} by
an arbitrary solution $h$ to equation \eqref{e4} in a ball $B$,
containing $D$, integrate by parts, take into account
the boundary conditions in \eqref{e2}, and get the relation
$$\int_D h(x)dx=0 \qquad \forall h\in \mathcal{N}.$$
Since $h\in \mathcal{N}$ implies $h(gx)\in \mathcal{N}$ forall $g$,
where $g$ is an arbitrary rotation in $\R^3$ about the origin $O$,
one obtains
\begin{equation}
\label{e8}
\int_D h(gx)dx=0 \qquad \forall h\in \mathcal{N},\,\, \forall g.
\end{equation}
Let $O\in D$,  and take an arbitrary straight
line $\ell$ passing through $O$ and directed along a
unit vector $\alpha$. Let
$g=g(\phi)$ be the rotation about
$\ell$ by an angle $\phi$ counterclockwise. Differentiate
\eqref{e8} with respect to $\phi$ and then set $\phi=0$, see \cite{R382}.
A similar argument has been used in \cite{R512}-\cite{R614}. The result 
is
\begin{equation}
\label{e9}
\int_D\nabla h(x)\cdot[\alpha, x]dx=0,
\end{equation}
where $[\alpha, x]$ is the cross product and $\cdot$ stands for the inner
product in $\R^3$. Equation  \eqref{e9} is invariant with respect to
translations because $\int_D\nabla h(x)dx=0$ $\forall h\in  \mathcal{N}$.
Indeed,  $h\in  \mathcal{N}$ implies $\nabla h \in  \mathcal{N}$.
Using the divergence theorem, the relation $\nabla \cdot [\alpha, x]=0$,
valid for any constant vector $\alpha$, and the arbitrariness of $\alpha$,
one derives from  \eqref{e9} the following relation
\begin{equation}
\label{e10}
\int_S h(s) [s, N]ds=0,\quad \forall h\in \mathcal{N},
\end{equation}
which is also invariant with respect to translations. Indeed, $N$ is
invariant under translations because $[s_p, s_q]$ is, and $\int_S
h(s) [a, N]ds=0$ for any constant vector $a$ because $\int_S h(s)
Nds=\int_D \nabla hdx=0$, as was pointed out above. Let us derive
from \eqref{e10} equation \eqref{e3}. We need the following result.

{\bf Lemma 1.}{\it  The orthogonal complement of the set $M$ of the
restrictions of all $h\in \mathcal{N}$ to $S$ is a finite-dimensional
space spanned by the functions $u_{jN}$, where $\{u_j\}_{j=1}^J$
is the basis of the eigenspace of the Dirichlet Laplacian in
$D$, corresponding to the eigenvalue $k^2$.}

{\bf Remark 1}. It follows from \eqref{e10} that $[s,N]$ is orthogonal
in $L^2(S)$ to the set $M$. Therefore, by Proposition 1,  each of
the three components of $[s,N]$ must be  linear combinations of the
functions $u_{jN},\,\,1\le j\le J$. In other words, equation
\eqref{e3} holds.

{\it Proof of Lemma 1.} Note that the  result of Lemma 1 is
equivalent to the assertion that the boundary value problem
\begin{equation}
\label{e11}
(\nabla^2+k^2)h=0\quad in \quad D, \quad h|_S=f
\end{equation}
is solvable if and only if
\begin{equation}
\label{e12}
 \int_Sfu_{jN}ds=0, \quad 1\le j\le J.
\end{equation}
where $u_j$, $1\le j\le J,$ is a basis of the solutions to problem
\eqref{e4}.

The {\it necessity}  of conditions \eqref{e12} is proved by the relation
$$0=\int_Du_j(\nabla^2+k^2)hdx=-\int_Sfu_{jN}ds,$$
where an integration by parts  and the boundary condition $u_j=0$ on
$S$ were used, and equation \eqref{e4} for $h$
was taken into account.

The {\it sufficiency} of conditions \eqref{e12} is proved as follows.
Denote by $H^{m}(D)$ the usual Sobolev spaces. Given
an $f\in H^{3/2}(S)$, construct
an arbitrary $F\in H^2(D)$, such that $F|_S=f|_S$, and define
$h:=w+F$, where
\begin{equation}
\label{e13}
(\nabla^2+k^2)w=-(\nabla^2+k^2)F \quad in \quad D, \quad w|_S=0.
\end{equation}
If such $w$ exists, then $h=w+F$ solves problem \eqref{e11}.
For the existence of $w$ it is necessary and sufficient that
$$\int_D(\nabla^2+k^2)Fu_jdx=0,\quad 1\le j\le J.$$
An integration by parts shows that
these conditions are equivalent to conditions \eqref{e12} because
$u_j$ solve problem \eqref{e4}. Thus, Lemma 1 is
proved. \hfill $\Box$

Equation \eqref{e10} says that $[s,N]$ is orthogonal to the
set $M$, that is, to the restrictions of all $h\in \mathcal{N}$ to $S$.

{\bf Lemma  2}. {\it The set $M$ is dense in $L^2(S)$ in the set
of
all $\psi\in H^{3/2}(S)$ for which the boundary problem \eqref{e13} is
solvable.}

This Lemma,  equation \eqref{e10}, and
Lemma 1 imply  \eqref{e3}.

{\it Proof of  Lemma 2.}
Assume the contrary. Then for some $f\in H^{3/2}(S)$, $f\neq 0$,  problem
\eqref{e11} is solvable and
$$\eta(y):=\int_S f(s)\psi(s,y)ds=0\quad \forall y\in D':=\R^3\setminus
D,$$ where $\psi(x,y):=\frac{e^{ik|x-y|}}{4\pi |x-y|}\in
\mathcal{N}$ for $y\in B'$, that is, outside a ball containing $D$.
The function  $\eta$ is a simple-layer potential which vanishes in
$B'$, and by the unique continuation property for solutions of the
homogeneous Helmholtz equation, $\eta=0$ everywhere in $D'$. Thus,
it vanishes on $S$. Therefore, $\eta$ solves problem \eqref{e4}. By
the jump relation for the normal derivative of $\eta$ across $S$,
one has $f=\eta_{N}$, where $\eta_{N}$ is the limiting value of the
normal derivative of $\eta$ on $S$ from inside $D$. If problem
\eqref{e10} is solvable, then, as we have proved, $f$ is orthogonal
to all functions $u_{jN}$. The function $\eta_{N}$ is a linear
combination of these functions. This and the relation $f=\eta_{N}$
prove that $f=0$. Consequently, we have proved the claimed density
of $M$ in the set of all $H^{3/2}(S)-$functions $f$ for which
problem \eqref{e11} is solvable.  Lemma 2 is proved. \hfill $\Box$

This completes the proof of Theorem 2. \hfill $\Box$

{\bf Theorem 3}. {\it Suppose {\it Assumptions A} hold and relation 
\eqref{e1} holds for some $f\neq 0$.  
Then problem \eqref{e2} has a solution. Conversely, if problem \eqref{e2} 
has a 
solution, then there exists $f\neq 0$ such that relation \eqref{e1} 
holds.}

{\it Proof of Theorem 3.} Write \eqref{e1} as
$$\int_{\R^3}f(gx+y)\chi (x)dx=0\qquad \forall y\in \R^3,\quad \forall
g\in G,$$ where $\chi (x)$ is the characteristic function of $D$.
Applying the Fourier transform and the convolution theorem one gets
\begin{equation}
\label{eX}
\tilde{f}(\xi)\overline{\tilde{\chi}(g^{-1}\xi)}=0,
\end{equation}
 where
$\tilde{f}$ and $\tilde{\chi}$ are the Fourier transforms of $f$ and
$\chi$, respectively, and the overbar stands for the complex
conjugate. The Fourier transform of $f$ is understood in the sense
of distributions. The Fourier transform $\tilde{\chi}$ is an entire
function of exponential type because function $\chi$ has support
$\overline{D}$, which is a bounded set. Moreover, $\tilde{\chi}$ is
a uniformly bounded function of $\xi\in \R^n$. The
product of the tempered distribution $\tilde{f}$ and the function
$\tilde{\chi}$ is a tempered distribution also, that is, and element of
$\mathcal{S}'$.

 Since $g^{-1}$
runs through all the rotations, one can replace $g^{-1}$ by $g$. It
follows from \eqref{eX} that
\begin{equation}
\label{eY}
supp\, \tilde{f}=\cup_{k} C_k,\,\, \text {where} \,\, C_k:=\{\xi:
\tilde{\chi}(\xi)=0\,\,\, \forall \xi: \xi^2-k^2=0\}.
\end{equation}
 In other
words, the support of the distribution $\tilde{f}$ is a subset of
the union of spherical surfaces of zeros of $\tilde{\chi}$, the
Fourier transform of the characteristic function of the bounded
domain $D$.
 Since $\tilde{\chi}(\xi)$ is an entire
function of exponential type, vanishing on an irreducible algebraic
variety $\xi^2-k^2=0$ in $\C^3$, one concludes, using the division
lemma,  that
\begin{equation}
\label{eZ}
\tilde{\chi}(\xi)=(\xi^2-k^2)\tilde{u}(\xi),
\end{equation}
 where $\tilde{u}$ is
an entire function of the same exponential type as $\tilde{\chi}$ (
see \cite{F}). Therefore, by the Paley-Wiener theorem, the
corresponding $u$ has compact support. Taking the inverse Fourier
transform of equation \eqref{eZ}, one gets:
\begin{equation}
\label{eW}
(-\nabla^2-k^2)u(x)=\chi(x)   \qquad in\quad \R^3, \quad u=0\quad
\text{if}\,\,
\quad |x|>R,
\end{equation}
 where $R>0$ is sufficiently large. By the elliptic
regularity results,  one concludes that  $u\in H^2_{loc}(\R^3)$.
Since $u$ solves the Helmholtz elliptic equation and vanishes near
infinity, that is, in the region  $|x|>R,$ the uniqueness of the
solution to the Cauchy problem to the equation \eqref{eW} and  the
path connectedness of the complement $D_1:=\bar{D}'$ of the closure
$\bar{D}$ of $D$ allow one to conclude that $u=0$ in $D_1$. The
connectedness and path connectedness of $D_1$ follow from our
Assumptions A and from the isotopy extension theorem (see \cite{H}). If
$u=0$ in $D_1$ and $u\in H^2_{loc}(\R^3)$, it follows from the
Sobolev embedding theorem that the boundary conditions (2) hold.
Since $\chi(x)=1$ in $D$, equation (2) holds. 
The converse statement in Theorem 3 has already been established
above, in our discussion of the equivalence of Conjecrues 1 and 2.
Theorem 3 is proved.
\hfill $\Box$.

\section{Relation to analyticity}

The classical Morera theorem in complex analysis says that if
$\int_Cf(z)dz=0$ for any closed polygon $C$ in a domain $D$ of the
complex plane, and if $f$ is continuously differentiable in $D$,
then $f$ is analytic in $D$. A simple proof is based on a version of
Green's formula:
$$0=\int_Cf(z)dz= 2i \int_{\Delta}\bar{\partial}f
dxdy.$$
Here $\Delta$ is the plane domain with the boundary $C$ and
$\bar{\partial}f:=\frac{f_x+if_y}{2}$. If
$\int_{\Delta}\bar{\partial}f dxdy=0$ for any polygon $\Delta$, then
one passes to the limit in the formula
$$\frac 1{|\Delta|}\int_{\Delta}\bar{\partial}f dxdy=0,$$
where $|\Delta|$ is the area of $\Delta$ and the limit is
taken as $diam \Delta\to 0$, so that $\Delta$ shrinks uniformly in 
directions to a point
$(x,y)\in \Delta$. Then for almost all points in $D$ one gets
$\bar{\partial}f=0$, and if $\bar{\partial}f$ is continuous, then
$\bar{\partial}f=0$ everywhere in $D$. This
implies that $f$ is analytic in $D$.

One may ask if the assumption that $f$ is continuously
differentiable can be replaced by a weaker assumption, and if the
set of polygons can be replaced by some other sets. The answer to
the first question is easy: if $f\in L^1_{loc}(D)$, then one
considers a mollified function $f_\epsilon(z):=\int_{\zeta:
|z-\zeta|\le \epsilon}\omega_\epsilon(z-\zeta)f(\zeta)dudv$, where
$\zeta=u+iv$ and $\omega_\epsilon(z)$ is the standard mollifying
kernel (\cite{Ho}, p.14). It is known that $f_\epsilon\to f$ in
$L^1(D)$ as $\epsilon \to 0$, and one can select a subsequence
$\epsilon_j\to 0$, such that $f_{\epsilon_j}\to f$ almost everywhere
in $D$. If $\int_Cf(z)dz=0$ for any closed polygon $C$, then
$\int_Cf_\epsilon(z)dz=0$ for any closed polygon $C$, and the above
argument, applied to the $C^1-$smooth $f_\epsilon$, leads to the 
conclusion that
$f_\epsilon$ is analytic in $D$ for all sufficiently small
$\epsilon$. Since a sequence $f_{\epsilon_j}$ of analytic functions
converges to $f$ in $L^1(D)$ and almost everywhere in $D$, one
concludes that $f$ is analytic in $D$. This follows from the
closedness of the differential operator $\bar{\partial}$. Namely,
one has $||f_{\epsilon_j}-f||_{L^1(D)}\to 0$ and
$\bar{\partial}f_{\epsilon_j}=0$, so
$||\bar{\partial}f_{\epsilon_j}-0||_{L^1(D)}\to 0$. Consequently,
$f$ belongs to the domain of the operator $\bar{\partial}$,  and
$\bar{\partial}f=0$ in $L^1(D)$. Therefore, $f$  is analytic in $D$.

The second question:
can one replace the set of polygons by other sets is less simple.
For example,  one cannot replace polygons by the set $\sigma(B)$,
where $B$ is a ball. Indeed, using the above argument one arrives
at the relation $\int_{\sigma(B)}\bar{\partial}fdxdy=0$, and this does not
imply that $\bar{\partial}f=0$, as the example on p. 2 shows.
 However, any domain $D$ which has $P$-property can be used in a
generalization of the Morera theorem. By $\partial D$ the boundary of
$D$ is denoted.

{\bf Theorem 4.} {\it Assume that $D$ has $P$-property, $f\in
L^1_{loc}(\R^2)$, and 
$$\int_{\partial \sigma(D)}fdz=0 \quad \forall \sigma\in G.$$ Then
$f$ is an entire function.}

{\it Proof.} By the argument given above it is sufficient to
prove this theorem assuming $f$ $C^1-$smooth. In this case
one has
\begin{equation}\label{eG}
\int_{\partial \sigma(D)}fdz=2i\int_{\sigma(D)}\bar{\partial}fdxdy=0,
\end{equation}
where $z=x+iy$.
Since $D$ has $P$-property, one concludes from the above equation that
$\bar{\partial}f=0$ in $\R^2$. This means that $f$ is an entire function.
\hfill $\Box$

Let $B_r$ denote a ball (disc if $n=2$) of radius $r$ centered at the 
origin, $f\in L^1_{loc}(\R^2)$, and $s_{j}$, $j=1,2,\ldots $, denote
positive zeros of the Bessel function $J_1(s)$. 

{\it From here to the end of 
Section 3 \it it is assumed that $n=2$}.

 In \cite{Z1} the following result is proved:

 {\it  If $\int_{\partial
\sigma(B_r)}fdz=0$ for $r=r_1$ and for $r=r_2$, and if $r_1/r_2$
does not belong to the set $s_{j}/s_{m}$ for any positive integers
$j$ and $m$, then $f$ is an entire function.}

This result is an immediate  consequence of Theorem 4. Indeed, it
follows from Example 1 with $n=2$ that if \eqref{e1} holds for 
$D=B_{r_1}$, then the support of $f$ belongs to the set
 $\mathcal{N}^{(1)}=\{\xi: |\xi|=s_{j}/r_1, \xi\in \R^2\,\, \text{for 
some 
positive integers}\,\, j$\}, see the proof of Theorem 3. 

Similarly, if \eqref{e1} holds for
$D=B_{r_2}$, then the support of $f$ belongs to the set
 $\mathcal{N}^{(2)}=\{\xi: |\xi|=s_{m}/r_2, \xi\in \R^2\,\, \text{ for 
some 
positive integers}\,\,
m$\}. If these two sets, $\mathcal{N}^{(1)}$ and 
$\mathcal{N}^{(2)}$, 
have 
empty intersection, then the support of $f$ is empty, so that $f=0$.
The role of $f$ will be played by  $\overline{\partial}f$ in what 
follows.
 
The set $\mathcal{N}^{(1)}$ does not intersect the set $\mathcal{N}^{(2)}$
if and only if $s_{j}/r_1\neq s_{m}/r_2$ for any
positive integers $j$ and $m$. This condition is equivalent to
the condition that $r_1/r_2$ does not belong to the set $s_{j}/s_{m}$
for any positive integers $j$ and $m$.
Under this condition one concludes that  
$\overline{\partial}f=0$, because the role of $f$ is played in our case
by $\overline{\partial}f$, see equation \eqref{eG}. 
Thus,  $f$ is an entire function. \hfill$\Box$ 

The proof of this result in \cite{Z1} is much longer and more complicated.

\section{ Another approach to Pompeiu problem and some remarks}

In Conjecture 3 we assume that
\begin{equation}\label{eE}
\int_D e^{ik\alpha \cdot x}dx=0 \qquad \forall \alpha\in S^2,\quad
k=const>0.
\end{equation}
In the derivation of this equation $\alpha$ can be an arbitrary
complex vector $z\in M$, where $M\subset \C^3$ is an algebraic
variety defined by the equation $z\cdot z=1$. Here $z\in \C^3$ and
$z\cdot w:=\sum_{j=1}^3 z_jw_j.$ Note that $w_j$ is used in the
definition of $z\cdot w$, and not its complex conjugate $\bar{w_j}$.
Let $z=a+ib$, where $a,b\in \R^3$. One checks easily that $z=a+ib\in
M$ if and only if $a\cdot b=0$ and $a^2-b^2=1$, where $a^2:=a\cdot
a$. Let $a=(\lambda^2+1)^{1/2}(e_1 \cos\theta +e_2 \sin \theta)$,
where $e_j$, $j=1, 2, 3,$ are unit vectors of a Cartesian basis in
$\R^3$, and $\theta\in [0,2\pi)$. Let $b=\lambda e_3$. Here
$\lambda\in \R$ is an arbitrary number. One can easily check that
$\pm a\pm ib\in M$.

Equation \eqref{eE} becomes
\begin{equation}\label{eE1}
\int_D e^{\pm \lambda k x_3+ik(\lambda^2+1)^{1/2}(x_1 \cos\theta +x_2
\sin\theta)}dx_1 dx_2 dx_3=0.
\end{equation}
Assume that $D\subset \R^3$ is a bounded domain  diffeomorphic to a ball.

{\bf Conjecture 6.} {\it Under these assumptions equation \eqref{eE1}
holds for all $\lambda\in \R^1$ and all $\theta\in [0,2\pi)$ if and only 
if $D$ is a ball.}

In order to prove Conjecture 6 it may help to assume additionally that $D$
is a convex centrally symmetric domain, but the author thinks that
Conjecture 6 is correct without additional assumptions.

If $D$ is a ball of radius $R$, then equation \eqref{eE1} holds if
$kR=s_{\frac 3 2, j}$, where $s_{\frac 3 2, j}$ are positive zeros
of the Bessel function $J_{3/2}$. This follows from the calculations
given in Example 1. Conjecture 6, if it is correct, says that
equation \eqref{eE1} cannot hold for any domain, satisfying the
stated assumptions, except for a ball. Asymptotic behavior, as
$\lambda\to \infty,$ of the integral in \eqref{eE1} may help to
verify Conjecture 6.

{\bf Remark 2.} It is proved in \cite{R363} and in \cite{R470} that
if $D_1\in P$ and $D_2$ is a "sufficiently close" to $D_1$, then
$D_2\in P$. This means that $P$-property is {\it stable} in some
sense. The  $\bar{P}$-property is {\it not stable}: small
perturbations of $D$ lead to domains the Fourier transform of the
characteristic function of which do not have a spherical surface of
zeros.

The notion of being "sufficiently close" is defined as follows.
The domain $D_2$ is $C^3$-smooth, strictly convex, its Gaussian curvature
is bounded from below by a positive constant, and $meas (D_{12}\setminus
D^{12})$ is sufficiently small. Here $D^{12}:=D_1\cap D_2$ and
$D_{12}:=D_1\cup D_2$.

{\bf Remark 3.} One can prove (see, for example, \cite{R470}, p.412) 
that if $D\subset \R^n$
is a bounded strictly convex domain with a smooth boundary, and
$\tilde{\chi}_D(t_m\alpha)=0$ for all $\alpha\in S^{n-1}$ and a sequence
$t_m\to +\infty$, then $D$ is a ball.

{\bf Remark 4.} It is easy to give examples of the domains $D$ which
have $P$-property: any polygon has this property since it does not
have a real-analytic boundary. An ellipsoid in $\R^n$, $n\ge 2$, has
$P$-property unless it is a ball. This is easy to check by
calculating the Fourier transform of the characteristic function of
an ellipsoid and checking that this Fourier transform does not have
a spherical set of zeros. If the equation of the ellipsoid $D$ is
$\sum_{j=1}^n \frac{x_j^2}{a_j^2}=1$, then this Fourier transform is
$$\tilde{\chi}_D=(2\pi)^{n/2}(\sum_{j=1}^n \xi_j^2a_j^2)^{-\frac {n}{2}}
J_{n/2}\big((\sum_{j=1}^n \xi_j^2a_j^2)^{1/2}\big).$$

This is calculated by making the change of variables $x_j'=x_j/a_j$,
which transforms the ellipsoid into the ball of radius $1$ and the
$\xi_j$ variable of the Fourier transform into $\xi_j a_j$, so that
then the formula for the Fourier transform of the characteristic
function of the ball can be used. This formula is given in Example 1.

{\bf Remark 5.} One can construct $f\neq 0$ satisfying equation
\eqref{e1} for the domain $D$, which fails to have $P$-property, by
the following method, see \cite{R470}, p.406. Let $|\xi|=b$ be the 
spherical surface $S_b$
of zeros of the Fourier transform of the characteristic function of
$D$. Take any function $A(\xi)\in L^1(S_b)$ and define
$\tilde{f}(\xi)=A(\xi)\delta(|\xi|-b)$, where $\delta(|\xi|-b)$ is
the delta-function supported on the sphere $S_b$. Then the inverse
Fourier transform of $\tilde{f}$ is a function $f\neq 0$, which
satisfies equation \eqref{e1}. Since $\tilde{f}$ has compact
support, by the Paley-Wiener theorem the function $f$ is an entire
function of $x$. For example, if $n=3$ and $A(\xi)=1$, then
$$f(x)=\frac 1 {(2\pi)^3}\int_{|\xi|=b}e^{-i\xi \cdot x}d\xi=
(2 \pi)^{-\frac 3 2}b^2\frac{J_{1/2}(b|x|)}{\sqrt{b|x|}},$$
 where the known formula
$$e^{ik\alpha \cdot x}=\sum_{\ell=0}^\infty 4\pi i^{\ell}j_{\ell}(k|x|)
Y_{\ell}(\alpha)Y_{\ell}(x^0)$$ was used. In this formula $k>0$ is a
constant, $\alpha\in S^2$ is a unit vector, $S^2$ is the unit sphere
in $\R^3$, $Y_{\ell}$ are the normalized in $L^2(S^2)$ spherical
harmonics, $x^0:=x/|x|$, and
$j_{\ell}(r):=\sqrt{\frac{\pi}{2r}}J_{\ell +\frac 1 2}(r)$.
$$ $$
{\bf Acknowledgments.} The author thanks Professors R. Burckel and  
C.N. Moore for reading the paper and comments.


Department of Mathematics, Kansas State University, Manhattan, KS
66506-2602,

email:     ramm@math.ksu.edu

\end{document}